\documentclass[12pt]{amsart}%
\usepackage{amsfonts}
\usepackage{latexsym}
\usepackage{amsmath}
\usepackage{amsfonts}
\usepackage{amstext}
\usepackage{amsxtra}
\usepackage{enumerate}
\usepackage{amsthm}
\usepackage{a4}
\usepackage{graphicx}
\usepackage{amssymb}
\usepackage{mathrsfs}%
\providecommand{\U}[1]{\protect\rule{.1in}{.1in}}

\newcommand{\K}{{\mathbb K}}

\newcommand{\norm}[1]{\|#1\|}

\theoremstyle{definition}
\newtheorem{definition}{Definition}[section]
\newtheorem{remark}[definition]{Remark}

\theoremstyle{plain}
\newtheorem{lemma}[definition]{Lemma}
\newtheorem{theorem}[definition]{Theorem}
\newtheorem{corollary}[definition]{Corollary}
\newtheorem{proposition}[definition]{Proposition}
\begin{document}
\title[New inclusion and coincidence theorems]{New inclusion and coincidence theorems for summing multilinear mappings}
\author{G. Botelho}
\address[Geraldo Botelho]{Faculdade de Matem\'atica, Universidade Federal de
Uberl\^andia, 38.400-902 - Uberl\^andia, Brazil}
\email{botelho@ufu.br}
\author{C. Michels}
\address[Carsten Michels]{Institut f\"ur Mathematik, Universit\"{a}t Oldenburg,
Postfach 2503, D-26111 Oldenburg, Germany}
\email{michels@mathematik.uni-oldenburg.de}
\author{D. Pellegrino}
\address[Daniel Pellegrino]{Depto de Matem\'{a}tica, UFPB, J. Pessoa, 58051-900, PB, Brazil}
\email{dmpellegrino@gmail.com}
\thanks{The third named author is supported by CNPq Grant 308084/2006-3.}

\begin{abstract}
In this paper we obtain new inclusion and coincidence theorems for absolutely
or multiple summing multilinear mappings. In particular, we derive optimal coincidence theorems of Bohnenblust-Hille type for multilinear forms on $K$-convex Banach spaces of cotype~$2$.
\end{abstract}
\maketitle

\section{Introduction and background}

For linear operators it is well-known that if $1\leq p\leq q<\infty$ then
every absolutely $p$-summing operator is absolutely $q$-summing; this type of
result is called \textquotedblleft inclusion theorem\textquotedblright. For
multilinear mappings the situation is different and there is no similar
result, in general. Such results for multilinear mappings have been
investigated by several authors in the recent years (see \cite{BBJP-PAMS, JMP,
inclusion, popa}) and the present paper presents new contributions in this direction.

Recently, in \cite{BBPR, BBJP-PAMS, JMP}, complex interpolation arguments were
used in order to obtain inclusion and coincidence theorems for spaces of
absolutely summing and multiple summing mappings involving spaces of cotype
$2;$ the interpolation results were based on the paper \cite{DM}. In this
paper, among other results, we investigate similar results for the cases of
spaces with cotype greater than $2$ as well as for $\mathcal{L}_{\infty}%
$-spaces. We also show that in some situations our results are optimal.

The roots of this line of investigation for multilinear mappings (and of the
search for coincidence results) can be traced back to Littlewood's celebrated
$4/3$ theorem for bilinear forms \cite{littlewood} and its following generalization due to Bohnenblust and Hille
\cite{BH}:

If $A\colon c_{0}\times\cdots\times c_{0}\rightarrow\mathbb{K}$ is a
continuous $n$-linear form, then there is a constant $C_{n}$ (depending only
on $n$) such that%
\begin{equation}
\left(
{\displaystyle\sum\limits_{i_{1},...,i_{n}=1}^{\infty}}
\left\vert A(e_{i_{1}},...,e_{i_{n}})\right\vert ^{\frac{2n}{n+1}}\right)
^{\frac{n+1}{2n}}\leq C_{n}\left\Vert A\right\Vert . \label{ttt}%
\end{equation}

The case $n=2$ recovers the classical Littlewood $4/3$ theorem.

In order to realize that Bohnenblust-Hille's theorem is in fact a predecessor
of today's coincidence theorems for multilinear mappings, it is worth
mentioning that a reformulation of (\ref{ttt}), due originally to
Per\'{e}z-Garc\'{\i}a, asserts that for every Banach spaces $E_{1}%
,\ldots,E_{n},$ every continuous $n$-linear form on $E_{1}\times\cdots\times
E_{n}$ is multiple $(\frac{2n}{n+1};1,\ldots,1)$-summing. The case
$E_{1}=\cdots=E_{n}=c_{0}$ recovers (\ref{ttt}) as a particular case. Using
the notation introduced below, this result can be stated as:
\begin{equation}
\label{bohnenblusthille}{\mathcal{L}}(E_{1},\ldots,E_{n};\mathbb{K}%
)={\mathcal{L}}_{ms(\frac{2n}{n+1};1,\ldots,1)}(E_{1},\ldots,E_{n}%
;\mathbb{K}),
\end{equation}
which makes clear what we mean by a coincidence theorem. Recently, the authors
in \cite{DS} have provided a unification and extension of several results related to the
original Bohnenblust-Hille result, in particular of a vector valued variant from \cite{bpgv}.

Our paper is organized as follows: After fixing some notation in section~\ref{notation}, we present our abstract  approach to inclusion and coincidence results for absolutely and multiple summing multilinear operators    using complex interpolation theory in section~\ref{absolutely} and section~\ref{multiple}, respectively. In section~\ref{bohisect}, we focus on Bohnenblust-Hille type results for multiple summing multilinear operators defined on spaces of finite cotype -- note that so far, most results of Bohnenblust-Hille type have been dealing with multilinear operators defined on $c_0$-spaces. A short appendix provides a clarification of some complexification arguments used throughout the paper.    

\section{Notation}
\label{notation}

Let $\mathbb{N}$ denote the set of natural numbers, $E,E_{1},\ldots,E_{n},F$
denote Banach spaces over $\mathbb{K}=\mathbb{R}$ or $\mathbb{C}$. For the
notions of cotype~$q \ge2$ and $\mathcal{L}_{p}$-space we refer to \cite{DJT}.
For $p\geq1,$ by $l_{p}(E)$ we mean the spaces of absolutely $p$-summable
sequences in $E;$ we represent by $l_{p}^{w}(E)$ the linear space of the
sequences $(x_{j})_{j=1}^{\infty}$ in $E$ such that $(\varphi(x_{j}%
))_{j=1}^{\infty}\in l_{p}$ for every continuous linear functional
$\varphi:E\rightarrow\mathbb{K}$. The map
\[
\Vert(x_{j})_{j=1}^{\infty}\Vert_{w,p}:=\sup_{\varphi\in B_{E^{\prime}}}%
\Vert(\varphi(x_{j}))_{j=1}^{\infty}\Vert_{p}%
\]
defines a norm in $l_{p}^{w}(E)$. When the sequences are finite (with $m$
terms) we write $l_{p}^{m}$ and $l_{p,w}^{m}$ instead of $l_{p}$ and
$l_{p}^{w},$ respectively. When $p=\infty$ we define%
\[
\Vert(x_{j})_{j=1}^{\infty}\Vert_{w,\infty}:=\sup\Vert x_{j}\Vert,
\]
i.e.,%
\[
l_{\infty}^{w}(E)=l_{\infty}(E).
\]
By $\mathcal{L}(E;F)$ we denote the Banach space of all bounded linear
operators between the Banach spaces $E$ and $F$, and by $\mathcal{L}%
_{as(p;q)}(E;F)$ the class of all absolutely $(p;q)$-summing linear operators
($1\leq q\leq p<\infty$), endowed with the usual norm $\left\Vert .\right\Vert
_{as(q;p)}.$ The space of all continuous $n$-linear mappings $A\colon
E_{1}\times\cdots\times E_{n}\rightarrow F,$ with the $\sup$ norm, will be
denoted by $\mathcal{L}(E_{1},\dots,E_{n};F)$. If $E_{1}=\cdots=E_{n}=E$ we
write $\mathcal{L}(^{n}E;F)$.

From now on, if $1\leq q<\infty,$ the symbol $q^{\ast}$ represents the
conjugate of $q$. It will be convenient to adopt that%
\[
\frac{q}{\infty}=0
\]
for any $q>0.$

\section{Absolutely summing multilinear operators}
\label{absolutely}

\begin{definition}
\textrm{Let $1\leq p_{1},\ldots,p_{n},q\leq\infty$ such that $1/q\leq1/{p_{1}%
}+\ldots+1/{p_{n}}$. An $n$-linear mapping $T\in\mathcal{L}(E_{1}%
,...,E_{n};F)$ is absolutely $(q;p_{1},...,p_{n})$-summing if there exists
$C\geq0$ such that
\begin{equation}
\left(  \sum\limits_{j=1}^{m}\parallel T(x_{j}^{(1)},...,x_{j}^{(n)}%
)\parallel^{q}\right)  ^{1/q}\leq C%
{\textstyle\prod\limits_{r=1}^{n}}
\left\Vert (x_{j}^{(r)})_{j=1}^{m}\right\Vert _{w,p_{r}},
\end{equation}
for every $m\in\mathbb{N}$ and $x_{j}^{(r)}\in E_{r}$, $j=1,...,m$ and
$r=1,...,n.$ For $q=\infty$, the left-handside has to be modified as usual,
taking the supremum over all $\parallel T(x_{j}^{(1)},...,x_{j}^{(n)}%
)\parallel.$}
\end{definition}

In this case we write $T\in\mathcal{L}_{as(q;p_{1},...,p_{n})}(E_{1}%
,...,E_{n};F)$, and $\pi_{as(q;p_{1},...,p_{n})}(T)$ denotes the infimum over
all $C$ as in the above. If $p_{1}=\cdots=p_{n}=p$, we write $\mathcal{L}%
_{as(q;p)}$ instead of $\mathcal{L}_{as(q;p,...,p)}.$ If $q=p_{1}=\cdots
=p_{n}$, we write $\mathcal{L}_{as,q}$ instead of $\mathcal{L}_{as(q;q,...,q)}%
$. The case $q=\infty$ clearly does not define anything new, but it will be
very helpful for interpolation purposes, since%
\[
\lbrack l_{\infty}^{m}(F),l_{p}^{m}(F)]_{\theta}=l_{r}^{m}(F)
\]
for%
\[
\frac{1}{r}=\frac{1-\theta}{p}%
\]
with isomorphism constant independent of $m.$ The next simple lemma will be
used several times along this paper:

\begin{lemma}
\label{aslemma} $\mathcal{L}_{as(\infty;p_{1},\ldots,p_{n})}(E_{1}%
,\ldots,E_{n};F)=\mathcal{L}(E_{1},\ldots,E_{n};F)$ with equal norms,
regardless of the choice of $1\leq p_{1},\ldots,p_{n}\leq\infty$.
\end{lemma}

Historically, the first coincidence result for absolutely summing multilinear
mappings is Defant-Voigt Theorem (see \cite[Theorem 3.10]{alencarmatos}), which we state below:

\begin{theorem}
[Defant-Voigt Theorem]\label{dvtheorem} For every Banach spaces $E_{1}%
,...,E_{n}$,
\[
\mathcal{L}_{as,1}(E_{1},\ldots,E_{n};\mathbb{K})=\mathcal{L}(E_{1}%
,\ldots,E_{n};\mathbb{K}).
\]
\end{theorem}

A first general inclusion formula can be found in \cite{matos}:

\begin{proposition}
\label{asinclusionmatos}
Let the indices $p \le q$ and $p_i \le q_i$, $i=1, \ldots,n$ be such that $0 \le \frac{1}{p_1}+ \ldots \frac{1}{p_n} - \frac{1}{p} \le \frac{1}{q_1}+ \ldots \frac{1}{q_n} -\frac{1}{q}$. Then $\mathcal{L}_{as(p;p_1, \ldots,p_n)} \subseteq \mathcal{L}_{as(q;q_1, \ldots,q_n)}$.
\end{proposition}

\subsection{Sandwich-type results}

\begin{proposition}
\label{assandwichmixed} Let $1\leq p\leq r\leq q\leq\infty$ and $1\leq
p_{1},\ldots,p_{n},q_{1},\ldots,q_{n},r_{1},\ldots,r_{n}\leq\infty$ such that
$1/t\leq1/{t_{1}}+\ldots+1/{t_{n}}$ for $t\in\{p,q,r\}$, and $T\in
\mathcal{L}(E_{1},\ldots,E_{n};F)$. Then $T\in\mathcal{L}_{as(p;p_{1}%
,\ldots,p_{n})}\cap\mathcal{L}_{as(q;q_{1},\ldots,q_{n})}$ implies
$T\in\mathcal{L}_{as(r;r_{1},\ldots,r_{n})}$, provided that there exists
$0<\theta<1$ such that $1/r=(1-\theta)/p+\theta/q$, and for all $i=1,\ldots,n$
it holds $1/{r_{i}}=(1-\theta)/{p_{i}}+\theta/{q_{i}}$ and one of the
following conditions is satisfied:

\begin{enumerate}
[(i)]

\item $E_{i}$ is an $\mathcal{L}_{\infty}$-space;

\item $E_{i}$ is of cotype~$2$ and $1 \le p_{i},q_{i} \le2$;

\item $E_{i}$ is of finite cotype~$s_{i}>2$ and $1\leq p_{i},q_{i}<s_{i}%
^{\ast}$;

\item $p_{i}=q_{i}=r_{i}$.
\end{enumerate}
\end{proposition}

\begin{proof}
We prove the complex case, and the real case then follows by complexification
as described in the appendix. Let $0<\theta<1$ be so that%
\begin{equation}
\frac{1}{r}=\frac{1-\theta}{p}+\frac{\theta}{q}. \label{hyt}%
\end{equation}
By assumption the map $T$ generates bounded operators%
\[
\left\{
\begin{array}
[c]{c}%
\widehat{T_{p}}:l_{p_{1},w}^{m}(E_{1})\times\cdots\times l_{p_{n},w}^{m}%
(E_{n})\rightarrow l_{p}^{m}(F)\\
\widehat{T_{q}}:l_{q_{1},w}^{m}(E_{1})\times\cdots\times l_{q_{n},w}^{m}%
(E_{n})\rightarrow l_{q}^{m}(F)
\end{array}
\right.
\]
Applying the complex interpolation method to these $n$-linear operators we get
a linear operator%
\[
\widehat{T_{(\theta)}}:\left[  l_{p_{1},w}^{m}(E_{1}),l_{q_{1},w}^{m}%
(E_{1})\right]  _{\theta}\times\cdots\times\left[  l_{p_{n},w}^{m}%
(E_{n}),l_{q_{n},w}^{m}(E_{n})\right]  _{\theta}\rightarrow\left[  l_{p}%
^{m}(F),l_{q}^{m}(F)\right]  _{\theta}%
\]
with%
\[
\left\Vert \widehat{T_{(\theta)}}\right\Vert \leq\left\Vert \widehat{T_{p}%
}\right\Vert ^{1-\theta}\left\Vert \widehat{T_{q}}\right\Vert ^{\theta}.
\]
This operator satisfies%
\[
\widehat{T_{(\theta)}}\left(  (x_{j}^{(1)})_{j=1}^{m},...,(x_{j}^{(n)}%
)_{j=1}^{m}\right)  =\left(  T(x_{j}^{(1)},...,x_{j}^{(n)})\right)  _{j=1}^{m}%
\]
for all sequences $(x_{j}^{(i)})_{j=1}^{m}$ in $\left[  l_{p_{i},w}^{m}%
(E_{i}),l_{q_{i},w}^{m}(E_{i})\right]  _{\theta},$ $1\leq i\leq n$.

By \cite[Theorem 5.1.2]{Bergh} we have $\left[  l_{p}^{m}(F),l_{q}%
(F)^{m}\right]  _{\theta}=l_{r}^{m}(F)$ isometrically with $r$ as in
(\ref{hyt}). Using the natural isometric identification $l_{p_{i},w}^{m}%
(E_{i})=l_{p_{i}}^{m}\otimes_{\varepsilon}E_{i}$, we will now see that for all
$i=1,\ldots,n$
\begin{equation}
l_{r_{i},w}^{m}(E_{i})=\left[  l_{p_{i},w}^{m}(E_{i}),l_{q_{i},w}^{m}%
(E_{i})\right]  _{\theta} \label{tensorinterpolation}%
\end{equation}
with isomorphism constant not depending on $m$. With this, we can identify the
operator $\widehat{T_{(\theta)}}$ with the map%
\[
\widehat{T_{r}}:l_{r_{1},w}^{m}(E_{1})\times\cdots\times l_{r_{n},w}^{m}%
(E_{n})\rightarrow l_{r}^{m}(F),
\]
and this gives us $T\in\mathcal{L}_{as(r;r_{1},\ldots,r_{n})}(E_{1}%
,\ldots,E_{n};F)$.

Coming back to \eqref{tensorinterpolation}, this is clear if condition~(iv) is
fulfilled. The apropriate statement under the assumptions in (ii) can be found
in \cite[Theorem]{DM}, for (iii) in \cite[Theorem~1]{Michels-Belg}. To see
that it holds under the assumptions in (i), first localize and then simply use
the fact that $l_{p_{i},w}^{m}(l_{\infty}^{k})=l_{p_{i}}^{m}\otimes
_{\varepsilon}l_{\infty}^{k}=l_{\infty}^{k}(\ell_{p_{i}}^{m})$ and complex
interpolation of vector-valued $l_{\infty}^{k}$'s. More precisely, we have%
\[
\lbrack l_{p_{i},w}^{m}(l_{\infty}^{k}),l_{q_{i},w}^{m}(l_{\infty}%
^{k})]_{\theta}=[l_{\infty}^{k}(\ell_{p_{i}}^{m}),l_{\infty}^{k}(\ell_{q_{i}%
}^{m})]_{\theta}=l_{\infty}^{k}(\ell_{r_{i}}^{m})=l_{r_{i},w}^{m}(l_{\infty
}^{k})
\]
for%
\[
\frac{1}{r_{i}}=\frac{1-\theta}{p_{i}}+\frac{\theta}{q_{i}},
\]
with isomorphism constant independent of $m$ and $k$.
\end{proof}

For $p=p_{1}=\ldots=p_{n}$, $q=q_{1}=\ldots=q_{n}$ and $r=r_{1}=\ldots=r_{n}$,
Theorem \ref{assandwichmixed} gives the following:

\begin{corollary}
\label{assandwich} Let $1\leq p<r<q\leq\infty$ and $T\in\mathcal{L}%
(E_{1},\ldots,E_{n};F)$. Then $T\in\mathcal{L}_{as,p}\cap\mathcal{L}_{as,q}$
implies $T\in\mathcal{L}_{as,r}$ in each of the following cases:

\begin{enumerate}
[(i)]

\item $E_{1},\ldots,E_{n}$ are all $\mathcal{L}_{\infty}$-spaces;

\item $1 \le p \le r \le q \le2$, and each $E_{i}$ is either an $\mathcal{L}%
_{\infty}$-space or of cotype~$2$, $i=1, \ldots,n$;

\item $1\leq p\leq r\leq q<s^{\ast}<2$, and each $E_{i}$ is either an
$\mathcal{L}_{\infty}$-space or of finite cotype~$s>2$, $i=1,\ldots,n$.
\end{enumerate}

\begin{proof}
Choose $0<\theta<1$ such that $1/r=(1-\theta)/p+\theta/q$ and apply
Proposition~\ref{assandwichmixed}.
\end{proof}
\end{corollary}

\subsection{Inclusion theorems for operators on spaces with finite cotype}

The following result extends \cite[Theorem 3]{JMP} to spaces with finite
cotype $>2.$ Note that while in the linear case we have a directed oriented
inclusion, that is,
\[
r\leq q\Longrightarrow{\mathcal{L}}_{as,r}(E;F)\subseteq{\mathcal{L}}%
_{as,q}(E;F),
\]
we are about to show that in the multilinear case the inclusion sometimes
holds in the opposite direction. It is worth noting that now we need the
hypothesis $n\geq s$ in (ii) of the next theorem.

\begin{theorem}
\label{opposite} Let $1 \le r < q<\infty$. Then $\mathcal{L}_{as,q}%
(E_{1},\ldots,E_{n};F) \subseteq\mathcal{L}_{as,r}(E_{1},\ldots,E_{n};F)$ if
one of the following conditions holds:

\begin{enumerate}
[(i)]

\item $1\leq r\leq q\leq2$, and $E_{1},\ldots,E_{n}$ of cotype~$2$ and
$n\geq2$;

\item $1\leq r\leq q<s^{\ast}<2$, $E_{1},\ldots,E_{n}$ of finite cotype~$s>2$,
and $n\geq s$.
\end{enumerate}
\end{theorem}

\begin{proof}
(i) is already known by \cite[Theorem~3]{JMP}. (ii) From \cite[Theorem 2.5]%
{B}, we know that
\[
\mathcal{L}_{as,1}(E_{1}, \ldots, E_{n};F)=\mathcal{L}(E_{1}, \ldots,
E_{n};F)\text{ }%
\]
holds true for $n\geq s$ and all Banach spaces $F$ provided that all $E_{i}$
have cotype $s$, $i=1,\ldots,n$. Now apply Corollary~\ref{assandwich} with
$p=1$.
\end{proof}

\begin{remark}
\textrm{\label{reft5} In the linear case we have $\mathcal{L}_{as,p}%
(E;F)=\mathcal{L}_{as,q}(E;F)$ whenever $E$ has finite cotype $s>2$ and $1\leq
p\leq q<s^{\ast}$ (cf. \cite[Corollary 11.16]{DJT}). } More recently, it was
shown in \cite{popa} and \cite{bparxiv2} (independently) that a similar
statement holds for multiple summing operators (for the notation we refer to
the next section): $\mathcal{L}_{ms,p}(E_{1},\ldots,E_{n};F)=\mathcal{L}%
_{ms,q}(E_{1},\ldots,E_{n};F)$ whenever all $E_{i}$ have finite cotype $s>2$
and $1\leq p\leq q<s^{\ast}$. Later on in subsection~\ref{reverse}, we will give an alternative proof of this result.  
\end{remark}

Using Theorem~\ref{dvtheorem}, a similar argument gives the following stronger
result for the scalar-valued case (note that here, there is no need of the
hypothesis $n\geq s$):

\begin{theorem}
\label{asscalarinclusion} Let $1 \le r < q<\infty$. Then $\mathcal{L}%
_{as,q}(E_{1},\ldots,E_{n};\mathbb{K}) \subseteq\mathcal{L}_{as,r}%
(E_{1},\ldots,E_{n};\mathbb{K})$ if one of the following conditions holds:

\begin{enumerate}
[(i)]

\item $E_{1},\ldots,E_{n}$ all $\mathcal{L}_{\infty}$-spaces;

\item $1 \le r \le q \le2$, and each $E_{i}$ is either $\mathcal{L}_{\infty}%
$-space or of cotype~$2$, $i=1, \ldots,n$;

\item $1 \le r \le q <s^{\ast}<2$, and each $E_{i}$ is either $\mathcal{L}%
_{\infty}$-space or of cotype~$s>2$, $i=1, \ldots,n$.
\end{enumerate}
\end{theorem}

\subsection{Inclusions and coincidences for operators on $\mathcal{L}_{\infty
}$-spaces}

For multilinear operators on $\mathcal{L}_{\infty}$-spaces, surprisingly the
usual inclusion in directed order holds, without any further assumptions on
the indices involved:

\begin{theorem}
\label{asinclusion} Let $1 \le p \le r \le\infty$ and $E_{1},\ldots,E_{n}$ all
$\mathcal{L}_{\infty}$-spaces. Then $\mathcal{L}_{as,p}(E_{1},\ldots,E_{n};F)
\subseteq\mathcal{L}_{as,r}(E_{1},\ldots,E_{n};F)$.
\end{theorem}

\begin{proof}
This follows immediately from Lemma~\ref{aslemma} and
Corollary~\ref{assandwich} with $q=\infty$.
\end{proof}

\begin{corollary}
\label{asinfinitycotype2} Let $E_{1},\ldots,E_{n}$ be all $\mathcal{L}%
_{\infty}$-spaces and $F$ a space of cotype~$2$. Then $\mathcal{L}%
(E_{1},\ldots,E_{n};F)=\mathcal{L}_{as,r}(E_{1},\ldots,E_{n};F)$ for all
$2\leq r\leq\infty$.
\end{corollary}

\begin{proof}
Use
$\mathcal{L}(E_{1},\ldots,E_{n};F)=\mathcal{L}_{as,2}(E_{1},\ldots
,E_{n};F)$ (\cite[Theorem 2.10]{Pell-Irish}) and
Theorem~\ref{asinclusion}.
\end{proof}

\begin{remark}
The above corollary also follows from the same (formally stronger) statement
for multiple summing operators later on (Corollary~\ref{msinfinitycotype2}).
\end{remark}

For the scalar field ${\mathbb{K}}$, the above can be strenghtened:

\begin{corollary}
Let $E_{1},\ldots,E_{n}$ be all $\mathcal{L}_{\infty}$-spaces. Then
$\mathcal{L}(E_{1},\ldots,E_{n};{\mathbb{K}})=\mathcal{L}_{as,r}(E_{1}%
,\ldots,E_{n};{\mathbb{K}})$ for all $1 \le r < \infty$.
\end{corollary}

\begin{proof}
This immediately follows from Theorem~\ref{dvtheorem} and
Theorem~\ref{asinclusion}.
\end{proof}

Surprisingly, we can even prove the following much stronger result,
that generalizes \cite[Corollary 2.5]{PPGG}:

\begin{theorem}
Let $E_{1},\ldots,E_{n}$ be all $\mathcal{L}_{\infty}$-spaces. Then
$\mathcal{L}(E_{1},\ldots,E_{n};{\mathbb{K}})=\mathcal{L}_{as(r;2r)}%
(E_{1},\ldots,E_{n};{\mathbb{K}})$ for all $1\leq r\leq\infty$.
\end{theorem}

\begin{proof}
Let $T\in\mathcal{L}(E_{1},\ldots,E_{n};{\mathbb{K}})$. From
\cite[Corolario~3.36]{David} or \cite[Corollary 2.5]{PPGG} we have
$T\in\mathcal{L}_{as(1;2)}.$ Recall that from Lemma \ref{aslemma} we also have
$T\in\mathcal{L}_{as(\infty;\infty)};$ so we use
Proposition~\ref{assandwichmixed} for $p=1,p_{i}=2$, $q=q_{i}=\infty$ and
$\theta=1-1/r.$
\end{proof}

\begin{remark}
Note that for $n=2$, one may deduce the above from the inclusion result
in Proposition~\ref{asinclusionmatos}, but not for $n\geq3.$
\end{remark}

\section{Multiple summing multilinear operators}
\label{multiple}

\begin{definition}
\textrm{Let $1 \le p_{1}, \ldots, p_{n} \le q \le\infty$. An $n$-linear
mapping $T\in\mathcal{L}(E_{1},...,E_{n};F)$ is multiple $(q;p_{1},...,p_{n}%
)$-summing if there exists $C\geq0$ such that
\begin{equation}
\left(  \sum\limits_{j_{1},...,j_{n}=1}^{m}\parallel T(x_{j_{1}}%
^{(1)},...,x_{j_{n}}^{(n)})\parallel^{q}\right)  ^{1/q}\leq C%
{\textstyle\prod\limits_{r=1}^{n}}
\left\Vert (x_{j}^{(r)})_{j=1}^{m}\right\Vert _{w,p_{r}},
\end{equation}
for every $m\in\mathbb{N}$ and $x_{j}^{(r)}\in E_{r}$, $j=1,...,m$ and
$r=1,...,n.$ For $q=\infty$, the left-handside has to be modified as usual,
taking the supremum over all $\parallel T(x_{j_{1}}^{(1)},...,x_{j_{n}}%
^{(n)})\parallel.$}
\end{definition}

In this case we write $T\in\mathcal{L}_{ms(q;p_{1},...,p_{n})}(E_{1}%
,...,E_{n};F)$, and $\pi_{ms(q;p_{1},...,p_{n})}(T)$ denotes the infimum over
all $C$ as in the above. If $p_{1}=\cdots=p_{n}=p$, we write $\mathcal{L}%
_{ms(q;p)}$ instead of $\mathcal{L}_{ms(q;p,...,p)}.$ If $q=p_{1}=\cdots
=p_{n}$, we write $\mathcal{L}_{ms,q}$ instead of $\mathcal{L}_{ms(q;q,...,q)}%
.$ The case $q=\infty$ clearly does not define anything new, but it will be
very helpful for interpolation purposes:

\begin{lemma}
\label{mslemma} $\mathcal{L}_{ms(\infty;p_{1},\ldots,p_{n})}(E_{1}%
,\ldots,E_{n};F)=\mathcal{L}(E_{1},\ldots,E_{n};F)$ with equal norms,
regardless of the choice of $1\leq p_{1},\ldots,p_{n}\leq\infty$.
\end{lemma}

It is worth mentioning that the theory of multiple summing mappings is quite
different from the theory of absolutely summing multilinear mappings and each
concept, in general, needs different techniques. Just to mention an example,
it is well known that in general
\[
\mathcal{L}%
(E_{1},\ldots,E_{n};\mathbb{K})\neq \mathcal{L}_{ms,1}(E_{1},\ldots,E_{n};\mathbb{K}),
\]
and this behavior is different from what asserts the Defant-Voigt Theorem for
absolutely summing multilinear forms.

\subsection{A Sandwich-type result}

\begin{proposition}
\label{mssandwichmixed} Let $1 \le p \le r \le q \le\infty$ and $1 \le
p_{1},\ldots,p_{n}, q_{1}, \ldots, q_{n}, r_{1}, \ldots,r_{n} \le\infty$ such
that ${t_{i}} \le t$ for $t \in\{p,q,r\}$ and all $i=1, \ldots,n$, and $T
\in\mathcal{L}(E_{1},\ldots,E_{n};F)$. Then $T \in\mathcal{L}_{ms(p;p_{1},
\ldots, p_{n})} \cap\mathcal{L}_{ms(q;q_{1}, \ldots, q_{n})}$ implies $T
\in\mathcal{L}_{ms(r;r_{1}, \ldots, r_{n})}$, provided that there exists $0
\le\theta\le1$ such that $1/r=(1-\theta)/p+\theta/q$, and for all $i=1,
\ldots,n$ it holds $1/{r_{i}}=(1-\theta)/{p_{i}} + \theta/{q_{i}}$ and one of
the following conditions is satisfied:

\begin{enumerate}
[(i)]

\item $E_{i}$ is $\mathcal{L}_{\infty}$-space;

\item $E_{i}$ of cotype~$2$ and $1 \le p_{i},q_{i} \le2$;

\item $E_{i}$ of cotype~$s_{i}>2$ and $1 \le p_{i},q_{i}<s_{i}^{*}$;

\item $p_{i}=q_{i}$.
\end{enumerate}
\end{proposition}

\begin{proof}
The proof for the case ${\mathbb{K}}={\mathbb{C}}$ goes along the same lines
as the one for absolutely summing multilinear operators in
Proposition~\ref{assandwichmixed}; just note that the exponent in the
vector-valued range space of the operators involved is $m^{n}$ instead of $m$.
Then complexification as described in the appendix proves the case
${\mathbb{K}}={\mathbb{R}}$.
\end{proof}

\begin{remark}
With regard to the inclusion theorem due to P\'{e}rez-Garc\'{\i}a
\cite{inclusion} which states that $\mathcal{L}_{ms,p}\subseteq\mathcal{L}%
_{ms,q}$ if $1\leq p<q<2$, it is superfluous to state analogs of
Corollary~\ref{assandwich}~(ii) and (iii); the case of $\mathcal{L}_{\infty}%
$-spaces will be dealt with in Section 4.3.
\end{remark}

\subsection{Reverse inclusions for multiple summing mappings}
\label{reverse}

Analogs of Theorem~\ref{opposite} for multiple summing operators have been
given recently in \cite{popa} and (independently) in \cite{bparxiv2}. In this
section we present a quite simple approach for these results, based only in
the linear theory.

\begin{lemma}
Let $1\leq r\leq2.$ If $E_{1},...,E_{n}$ have cotype $2$, then%
\[
\mathcal{L}_{ms,r}(E_{1},\ldots,E_{n};F)\subseteq\mathcal{L}_{ms,1}(E_{1}%
,\ldots,E_{n};F).
\]

\end{lemma}

\begin{proof}
Suppose that $A\in\mathcal{L}_{ms,r}(E_{1},\ldots,E_{n};F)$. Let $(z_{j}%
^{(k)})_{j}\in l_{1}^{w}(E_{k})$, $k=1,...,n.$ Since $E_{k}$ has cotype $2$,
we have%
\begin{equation}
\mathcal{L}_{as,r^{\ast}}(c_{0};E_{k})=\mathcal{L}(c_{0};E_{k})\label{31oc}%
\end{equation}
for every $k=1,...,n$. We know, from \cite[Proposition 2.2]{DJT}, that there
exist continuous linear operators%
\[
u_{k}:c_{0}\rightarrow E_{k}%
\]
so that $u_{k}(e_{j})=z_{j}^{(k)}$ for every $j$. Since $\left(  e_{j}\right)
_{j}\in l_{1}^{w}(c_{0})$, it follows from (\ref{31oc}) and \cite[Lemma
2.23]{DJT} that
\[
(z_{j}^{(k)})_{j}=(a_{j}^{(k)}y_{j}^{(k)})_{j},
\]
with $(a_{j}^{(k)})_{j}\in l_{r^{\ast}}$ and $(y_{j}^{(k)})_{j}\in l_{r}%
^{w}(E_{k}).$ Then%
\begin{align*}
&
{\displaystyle\sum\limits_{j_{1},...,j_{n}=1}^{\infty}}
\left\Vert A\left(  z_{j_{1}}^{(1)},...,z_{j_{n}}^{(n)}\right)  \right\Vert \\
&  =%
{\displaystyle\sum\limits_{j_{1},...,j_{n}=1}^{\infty}}
\left\Vert a_{j_{1}}^{(1)}...a_{j_{n}}^{(n)}A\left(  y_{j_{1}}^{(1)}%
,...,y_{j_{n}}^{(n)}\right)  \right\Vert \\
&  \leq\left(
{\displaystyle\sum\limits_{j_{1},...,j_{n}=1}^{\infty}}
\left\vert a_{j_{1}}^{(1)}...a_{j_{n}}^{(n)}\right\vert ^{r^{\ast}}\right)
^{1/r^{\ast}}\left(
{\displaystyle\sum\limits_{j_{1},...,j_{n}=1}^{\infty}}
\left\Vert A\left(  y_{j_{1}}^{(1)},...,y_{j_{n}}^{(n)}\right)  \right\Vert
^{r}\right)  ^{1/r}<\infty.
\end{align*}
Hence $A\in\mathcal{L}_{ms,1}(E_{1},\ldots,E_{n};F).$
\end{proof}

\begin{theorem}
If $E_{1},...,E_{n}$ have cotype $2$, then%
\[
\mathcal{L}_{ms,p}(E_{1},\ldots,E_{n};F)=\mathcal{L}_{ms,r}(E_{1},\ldots
,E_{n};F)
\]
for every $1\leq p\leq r<2.$
\end{theorem}

\begin{proof}
Proof. The inclusion $\subseteq$ is due to David-P\'{e}rez-Garc\'{\i}a
\cite{inclusion}. This result combined with the previous lemma completes the proof.
\end{proof}

With similar arguments one can show that

\begin{theorem}
If $E_{1},...,E_{n}$ have cotype $q>2$, then%
\[
\mathcal{L}_{ms,s}(E_{1},\ldots,E_{n};F)=\mathcal{L}_{ms,p}(E_{1},\ldots
,E_{n};F)
\]
for every $1\leq s\leq p<q^{\ast}.$
\end{theorem}

\begin{remark}
Using the estimates for the norms in \cite[Lemma 2.23]{DJT} it is possible to
get estimates for the norms of the spaces of multiple summing mappings in the
above theorems.
\end{remark}

\subsection{Inclusions and coincidences for multilinear operators on
$\mathcal{L}_{\infty}$-spaces}

As for absolutely summing multilinear operators in Theorem~\ref{asinclusion},
one can obtain the following inclusion result for multiple summing operators
on $\mathcal{L}_{\infty}$-spaces:

\begin{theorem}
\label{msinclusion} Let $1\leq p\leq r\leq\infty$ and $E_{1},\ldots,E_{n}$ all
$\mathcal{L}_{\infty}$-spaces. Then $\mathcal{L}_{ms,p}(E_{1},\ldots
,E_{n};F)\subseteq\mathcal{L}_{ms,r}(E_{1},\ldots,E_{n};F)$.
\end{theorem}

\begin{proof}
Choose $0<\theta<1$ such that $1/r=(1-\theta)/p$ and apply
Proposition~\ref{mssandwichmixed} (i) with $q=q_{1}=\cdots=q_{n}=\infty$ and
$p_{1}=\cdots=p_{n}=p$.
\end{proof}

The following (formally) improves upon Corollary~\ref{asinfinitycotype2}.

\begin{corollary}
\label{msinfinitycotype2} Let $E_{1},\ldots,E_{n}$ all be $\mathcal{L}%
_{\infty}$-spaces, and $F$ of cotype~$2$. Then $\mathcal{L}(E_{1},\ldots
,E_{n};F)=\mathcal{L}_{ms,r}(E_{1},\ldots,E_{n};F)$ for all $2 \le r \le
\infty$.
\end{corollary}

\begin{proof}
Just recall that under our hypothesis we have $\mathcal{L}(E_{1},\ldots
,E_{n};F)=\mathcal{L}_{ms,2}(E_{1},\ldots,E_{n};F)$ (\cite[Theorem~3.1]%
{bpgv}), and then use Theorem \ref{msinclusion}.
\end{proof}

\section{Bohnenblust-Hille type results}
\label{bohisect}

\subsection{A vector valued result}

Although in this section we mainly consider multilinear forms, let us start
with a more general result. We will need the following lemma due to David
P\'{e}rez-Garc\'{\i}a \cite[Teorema 5.2]{David} and Marcela Souza
\cite[Teorema 1.7.3]{souza} (see also \cite[Theorem~3.2]{bpgv}):

\begin{lemma}
\label{DJT1} If $F$ has cotype $q\geq2$, then for any Banach spaces
$E_{1},\ldots,E_{n}$ we have
\[
\mathcal{L}(E_{1},\ldots,E_{n};F)=\mathcal{L}_{ms(q;1)}(E_{1},\ldots,E_{n};F).
\]

\end{lemma}

Also, we need a refinement of \cite[3.16]{pgvck} (the proof is essentially the
same, so we omit it).

\begin{lemma}
\label{cklemma} Let $1 \le p_{1}, \ldots,p_{n} \le q<\infty$, and let $E_{i}$
for $i \in\{1,\ldots,n\}$ be an $\mathcal{L}_{\infty}$-space and $p_{i}<q$.
Then a multilinear operator $T \in\mathcal{L}(E_{1}, \ldots,E_{n};F)$ is
multiple $(q;p_{1},\ldots,p_{n})$-summing if and only if it is multiple
$(q;p_{1}, \ldots, \tilde{p}_{i},\ldots,p_{n})$-summing for all $1 \le
\tilde{p}_{i}<q$.
\end{lemma}

\begin{proposition}
\label{bhcotype} Let $F$ be of finite cotype~$q\geq2$. Then for $r>q$,
\[
\mathcal{L}(E_{1},\ldots,E_{n};F)=\mathcal{L}_{ms(r;p_{1},\ldots,p_{n})}%
(E_{1},\ldots,E_{n};F)
\]
provided that for each $i=1,\ldots,n$ one of the following conditions holds:

\begin{enumerate}
[(i)]

\item $E_{i}$ is an $\mathcal{L}_{\infty}$-space and $1 \le p_{i}< r $;

\item $E_{i}$ is of cotype~$2$ and $\frac{1}{p_{i}}=\frac{q}{2r}+\frac{1}{2}$;

\item $E_{i}$ is of cotype~$s_{i}>2$ and $\frac{1}{p_{i}}>\frac{1}{s_{i}^{*}%
}+\frac{q}{s_{i} \, r}$;

\item $p_{i}=1$.
\end{enumerate}
\end{proposition}

\begin{proof}
If $T\in\mathcal{L}(E_{1},\ldots,E_{n};F)$, then from Lemma~\ref{DJT1} and
Lemma \ref{mslemma} we have
\[
T\in\mathcal{L}_{ms(q;1,\ldots,1)}(E_{1},\ldots,E_{n};F)\cap\mathcal{L}%
_{ms(\infty,q_{1}\ldots,q_{n})}(E_{1},\ldots,E_{n};F)
\]
with
\[
q_{i}%
\begin{cases}
=2 & \text{ if }E_{i}\text{ has cotype }2\text{ or }E_{i}\text{ is an
}\mathcal{L}_{\infty}-space,\\
<s_{i}^{\ast} & \text{ if }E_{i}\text{ has cotype }s_{i}>2.
\end{cases}
\]
Then use Proposition~\ref{mssandwichmixed} with $\theta=1-\frac{q}{r}$. If
$E_{i}$ is an $\mathcal{L}_{\infty}$-space, we apply Lemma~\ref{cklemma} to
improve the corresponding summability index.
\end{proof}

\subsection{Multilinear forms on spaces with finite cotype}

For $F={\mathbb{K}}$, we can do much better. This requires a little preparation.

\begin{proposition}
\label{w}Let $q\geq2$ and $n\geq1.$ If $E_{1},\ldots,E_{n}$ are Banach spaces
and
\[
\mathcal{L}(E_{1},\ldots,E_{n};\mathbb{K})=\mathcal{L}_{ms(q;q_{1}%
,\ldots,q_{n})}(E_{1},\ldots,E_{n};\mathbb{K}),
\]
then for any Banach space $E_{n+1}$ we have%
\[
\mathcal{L}(E_{1},\ldots,E_{n+1};\mathbb{K})=\mathcal{L}_{ms(q;q_{1}%
,\ldots,q_{n},1)}(E_{1},\ldots,E_{n+1};\mathbb{K}).
\]

\end{proposition}

\begin{proof}
Let
\[
(x_{i_{r}}^{(r)})_{i_{r}=1}^{m}\subset E_{r},\text{ }r=1,...,n\text{ and
}(y_{j})_{j=1}^{m}\subset E_{n+1}%
\]
be given. For sake of abbreviation we put
\[
\mathbf{x}_{\mathbf{i}}=(x_{i_{1}}^{(1)},\ldots,x_{i_{n}}^{(n)})\text{ for
}\mathbf{i}=(i_{1},\ldots,i_{n}).
\]

Define $S\in\mathcal{L}(E_{n+1};l_{q})$ by
\[
S(y)=(T(\mathbf{x}_{\mathbf{i}},y)_{\mathbf{i}\in\{1,...,m\}^{n}}%
,0,0,\ldots)\in l_{q}.
\]
Since $l_{q}$ has cotype $q$ we have $S\in\mathcal{L}_{as(q;1)}(E_{n+1}%
;l_{q})$, and%
\[
\Vert S\Vert_{as(q;1)}\leq c_{0}\Vert S\Vert
\]
for some $c_{0}>0$ (not depending on $m$). Further, we have (using the
hypothesis)
\begin{align*}
\Vert S\Vert &  =\sup\limits_{y\in B_{E_{n+1}}}\Vert S(y)\Vert_{q}%
=\sup\limits_{y\in B_{E_{n+1}}}\left(  \sum\limits_{i}|T(\mathbf{x}%
_{i},y)|^{q}\right)  ^{1/q}\\
&  \leq\sup\limits_{y\in B_{E_{n+1}}}\Vert T(\cdot,y)\Vert_{ms(q;q_{1}%
,..,q_{n})}\prod\limits_{r=1}^{n}\Vert(x_{j_{r}}^{(r)})_{j_{r}=1}^{m}%
\Vert_{w,q_{r}}\\
&  \leq\sup\limits_{y\in B_{E_{n+1}}}c_{1}\Vert T(\cdot,y)\Vert\prod
\limits_{r=1}^{n}\Vert(x_{j_{r}}^{(r)})_{j_{r}=1}^{m}\Vert_{w,q_{r}}\\
&  \leq c_{1}\Vert T\Vert\prod\limits_{r=1}^{n}\Vert(x_{j_{r}}^{(r)}%
)_{j_{r}=1}^{m}\Vert_{w,q_{r}}%
\end{align*}
for some $c_{1}>0$ (not depending on $m$). So we have
\begin{align*}
\left(  \sum\limits_{\mathbf{i}}\sum\limits_{j}|T(\mathbf{x}_{\mathbf{i}%
},y_{j})|^{q}\right)  ^{1/q} &  =\left(  \sum\limits_{j}\Vert Sy_{j}\Vert
^{q}\right)  ^{1/q}\\
&  \leq\Vert S\Vert_{as(q;1)}\Vert(y_{j})_{j=1}^{m}\Vert_{w,1}\\
&  \leq c_{0}\Vert S\Vert\Vert(y_{j})_{j=1}^{m}\Vert_{w,1}\\
&  \leq c_{0}c_{1}\Vert T\Vert\left(  \prod\limits_{r=1}^{n}\Vert(x_{j_{r}%
}^{(r)})_{j_{r}=1}^{m}\Vert_{w,q_{r}}\right)  \Vert(y_{j})_{j=1}^{m}%
\Vert_{w,1}%
\end{align*}
which completes the proof.
\end{proof}

The following result appears in \cite[Proposition 3.5]{bparxiv1}. For the sake
of completeness we present a proof:

\begin{corollary}
\label{gg}Let $q\geq2$ and $n\geq2.$ If $E_{1},\ldots,E_{n}$ are Banach
spaces, then
\[
\mathcal{L}(E_{1},\ldots,E_{n};\mathbb{K})=\mathcal{L}_{ms(q;1,\ldots
,1,r)}(E_{1},\ldots,E_{n};\mathbb{K})
\]
for every $1\leq r\leq q.$
\end{corollary}

\begin{proof}
It is obvious that we just need to consider the case $r=q$.

We know that%
\[
\mathcal{L}(E_{n};\mathbb{K})=\mathcal{L}_{as(q;q)}(E_{n};\mathbb{K}).
\]
So, from Proposition \ref{w} we have%
\[
\mathcal{L}(E_{n-1},E_{n};\mathbb{K})=\mathcal{L}_{ms(q;1,q)}(E_{n-1}%
,E_{n};\mathbb{K}).
\]
By repeating this procedure we get%
\[
\mathcal{L}(E_{1},\ldots,E_{n};\mathbb{K})=\mathcal{L}_{ms(q;1,\ldots
,1,q)}(E_{1},\ldots,E_{n};\mathbb{K}).
\]

\end{proof}

We can now state our variant of the Bohnenblust-Hille result. Recall from the
introduction that this, by a suitable reformulation, says that ${\mathcal{L}%
}(E_{1},\ldots,E_{n};\mathbb{K})={\mathcal{L}}_{ms(\frac{2n}{n+1};1,\ldots
,1)}(E_{1},\ldots,E_{n};\mathbb{K})$ for all Banach spaces $E_{1},\ldots
,E_{n}$. Clearly, for $n \rightarrow\infty$, the first index tends to $2$. Our
result below gives information under which conditions on the spaces and the
indices involved every multilinear form is multiple $(2;p_{1},\ldots,p_{n})$-summing.

\begin{theorem}
\label{bohnenblust2} Let $E_{1}, \ldots,E_{n}$ be Banach spaces with finite
cotype. Then
\[
\mathcal{L}(E_{1},\ldots,E_{n};{\mathbb{K}})=\mathcal{L}_{ms(2;p_{1}^{(n)},
\ldots, p_{n}^{(n)})}(E_{1},\ldots,E_{n};{\mathbb{K}}),
\]
where
\[
p_{i}^{(n)}=
\begin{cases}
\frac{2n}{2n-1} & \text{if $E_{i}$ is of cotype~$2$,}\\
\frac{n q_{i,0}}{(n-1)q_{i,0} +1} & \text{if $E_{i}$ is of cotype~$q_{i}>2$
and $1 \le q_{i,0}<q_{i}^{\ast}$.}%
\end{cases}
\]

\end{theorem}

\begin{proof}
We are going to prove by induction over $n$. The case $n=1$ is trivial.

Suppose now that the result is true for some $n$. Let us consider any
$(n+1)$-linear form $T\in\mathcal{L}(E_{1},\ldots,E_{n},E_{n+1};\mathbb{K})$.
From Proposition \ref{w} and our hypothesis on $E_{1},\ldots,E_{n}$ we know
that
\begin{equation}
T\in\mathcal{L}_{ms(2;p_{1}^{(n)},\ldots,p_{n}^{(n)},1)}(E_{1},\ldots
,E_{n},E_{n+1};\mathbb{K})
\end{equation}
with $p_{i}^{(n)}$ as in the formulation of the theorem. From Corollary
\ref{gg} we know that%
\[
T\in\mathcal{L}_{ms(2;1,\ldots,1,r_{n+1})}(E_{1},\ldots,E_{n+1};\mathbb{K}),
\]
where we choose
\[
r_{n+1}:=%
\begin{cases}
2 & \text{if $E_{n+1}$ is of cotype~$2$},\\
q_{n+1,0} & \text{if $E_{n+1}$ is of cotype $q_{n+1}>2$ and $q_{n+1,0}%
<q_{n+1}^{\ast}$}.
\end{cases}
\]
Now we use Proposition~\ref{mssandwichmixed} with $\theta=\frac{n}{n+1}.$ Note
that for $i=1,\ldots,n$ it is
\[
\frac{1}{p_{i}^{(n+1)}}:=\frac{\theta}{p_{i}^{(n)}}+\frac{1-\theta}{1}=%
\begin{cases}
\frac{1}{\frac{2(n+1)}{2(n+1)-1}} & \text{if $E_{i}$ is of cotype~$2$},\\
\frac{1}{\frac{(n+1)q_{i,0}}{nq_{i,0}+1}} & \text{if $E_{i}$ is of cotype
$q_{i}>2$},
\end{cases}
\]
and
\[
\frac{1}{p_{n+1}^{(n+1)}}:=\frac{\theta}{1}+\frac{1-\theta}{r_{n+1}}=%
\begin{cases}
\frac{1}{\frac{2(n+1)}{2(n+1)-1}} & \text{if $E_{n+1}$ is of cotype~$2$},\\
\frac{1}{\frac{(n+1)q_{n+1,0}}{nq_{n+1,0}+1}} & \text{if $E_{n+1}$ is of
cotype $q_{n+1}>2$}.
\end{cases}
\]
So we get
\[
T\in\mathcal{L}_{ms(2;p_{1}^{(n+1)},\ldots,p_{n+1}^{(n+1)})}(E_{1}%
,\ldots,E_{n+1};\mathbb{K}),
\]
and the case $n+1$ is done.
\end{proof}

\begin{corollary}
\label{bhcotype2} \label{1}Let $n\geq1$ and let $E_{1},\ldots,E_{n}$ be Banach
spaces of cotype $2$. Then
\[
\mathcal{L}(E_{1},\ldots,E_{n};\mathbb{K})=\mathcal{L}_{ms(2;\frac{2n}{2n-1}%
)}(E_{1},\ldots,E_{n};\mathbb{K}).
\]

\end{corollary}

\begin{remark}
\begin{enumerate}
[(a)]

\item Note that the Defant-Voigt theorem together with the inclusion formula
for absolutely $(q;p_{1}, \ldots, p_{n})$-summing operators in
Proposition~\ref{asinclusionmatos} implies that
\[
\mathcal{L}(E_{1},\ldots,E_{n};\mathbb{K})=\mathcal{L}_{as(2;\frac{2n}{2n-1}%
)}(E_{1},\ldots,E_{n};\mathbb{K})
\]
for all Banach spaces $E_{1}, \ldots, E_{n}$. Thus, one can view the above
corollary as a variant of the Defant-Voigt theorem for multiple summing
operators on cotype~$2$ spaces.

\item In \cite[Theorem 2.3]{BBJP-PAMS} it was shown that $\mathcal{L}%
(E_{1},\ldots,E_{n};{\mathbb{K}})=\mathcal{L}_{ms(2;q_{k})}(E_{1},\ldots
,E_{n};{\mathbb{K}})$ with $q_{k}=\frac{2^{k+1}}{2^{k+1}-1}$ and $k$ such that
$2^{k-1}<n\leq2^{k}$ if $E_{1},\ldots,E_{n}$ are all of cotype~$2$. In
particular, for $n=2^{k}$ it is $q_{k}=\frac{2n}{2n-1}$. In fact, our theorem
above now shows that this is valid for $n$ arbitrary and that our estimates
improve the previous from \cite[Theorem 2.3]{BBJP-PAMS}. Just to give an
example, if $n=3$, \cite[Theorem 2.3]{BBJP-PAMS} gives that%
\[
\mathcal{L}(E_{1},E_{2},E_{3};\mathbb{K})=\mathcal{L}_{ms(2;\frac{8}{7}%
)}(E_{1},E_{2},E_{3};\mathbb{K}).
\]
On the other hand our result gives%
\[
\mathcal{L}(E_{1},E_{2},E_{3};\mathbb{K})=\mathcal{L}_{ms(2;\frac{6}{5}%
)}(E_{1},E_{2},E_{3};\mathbb{K}).
\]

\end{enumerate}
\end{remark}

In the case that all spaces involved have cotype~$q>2$, we get the following analog:

\begin{corollary}
\label{q}Let $n\geq1$ and let $E_{1},\ldots,E_{n}$ be Banach spaces of cotype
$q>2$. Then
\[
\mathcal{L}(E_{1},\ldots,E_{n};\mathbb{K})=\mathcal{L}_{ms(2;\frac{qn}%
{qn-1}-\varepsilon)}(E_{1},\ldots,E_{n};\mathbb{K})
\]
for all $\varepsilon>0$ sufficiently small.
\end{corollary}

\begin{proof}
By Theorem~\ref{bohnenblust2} one gets that $\mathcal{L}(E_{1},\ldots
,E_{n};\mathbb{K})=\mathcal{L}_{ms(2;\frac{nq^{\ast}}{(n-1)q^{\ast}%
+1}-\varepsilon)}(E_{1},\ldots,E_{n};\mathbb{K})$ for all $\varepsilon>0$
sufficiently small. An elementary calculation now shows that $\frac{nq^{\ast}%
}{(n-1)q^{\ast}+1}=\frac{qn}{qn-1}$.
\end{proof}

Further interpolation with the original Bohnenblust-Hille result and
Lemma~\ref{mslemma}, respectively, gives us the following more general statement:

\begin{corollary}
\label{bohnenblust3} Let $E_{1},\ldots,E_{n}$ be Banach spaces with finite
cotype, and let $\frac{2n}{n+1}\leq r<\infty$. Then
\[
\mathcal{L}(E_{1},\ldots,E_{n};{\mathbb{K}})=\mathcal{L}_{ms(r;r_{1}%
^{(n)},\ldots,r_{n}^{(n)})}(E_{1},\ldots,E_{n};{\mathbb{K}}),
\]
where for $i=1,\ldots,n$
\[
\frac{1}{r_{i}^{(n)}}=%
\begin{cases}
\frac{1}{2}+\frac{1}{r}-\frac{1}{\max(r,2)n} & \text{if $E_{i}$ is of
cotype~$2$,}\\
\frac{1}{q_{i,0}}+\frac{2}{r{q_{i,0}}^{\ast}}-\frac{2}{\max(r,2){q_{i,0}%
}^{\ast}n} & \text{if $E_{i}$ is of cotype~$q_{i}>2$ and $1\leq q_{i,0}%
<q_{i}^{\ast}$.}%
\end{cases}
\]

\end{corollary}

\begin{proof}
(i) The case $\frac{2n}{n+1} \le r \le2$:

By \eqref{bohnenblusthille} we know that
\[
{\mathcal{L}}(E_{1},\ldots,E_{n};\mathbb{K})={\mathcal{L}}_{ms(\frac{2n}%
{n+1};1,\ldots,1)}(E_{1},\ldots,E_{n};\mathbb{K}),
\]
and by Theorem~\ref{bohnenblust2}
\[
\mathcal{L}(E_{1},\ldots,E_{n};{\mathbb{K}})=\mathcal{L}_{ms(2;p_{1}%
^{(n)},\ldots,p_{n}^{(n)})}(E_{1},\ldots,E_{n};{\mathbb{K}}),
\]
where $p_{i}^{(n)}$, $i=1,\ldots,n$ are as in the statement of the theorem.
Now let $\frac{2n}{n+1}<r<2$. Then Proposition~\ref{mssandwichmixed} with
$\theta=(n+1)-\frac{2n}{r}$ gives
\[
\mathcal{L}(E_{1},\ldots,E_{n};{\mathbb{K}})=\mathcal{L}_{ms(r;r_{1}%
^{(n)},\ldots,r_{n}^{(n)})}(E_{1},\ldots,E_{n};{\mathbb{K}}),
\]
where for $i=1,\ldots,n$
\[
\frac{1}{r_{i}^{(n)}}=%
\begin{cases}
\frac{1-\theta}{1}+\frac{\theta(2n-1)}{2n}=\frac{1}{2}+\frac{1}{r}-\frac
{1}{2n} & \text{if $E_{i}$ is of cotype~$2$,}\\
\frac{1-\theta}{1}+\frac{\theta({q_{i,0}}^{\ast}n-1)}{{q_{i,0}}^{\ast}n}%
=\frac{1}{q_{i,0}}+\frac{2}{r{q_{i,0}}^{\ast}}-\frac{1}{{q_{i,0}}^{\ast}n} &
\text{if $E_{i}$ is of cotype~$q_{i}>2$}\\
& \text{and $1\leq q_{i,0}<q_{i}^{\ast}$.}%
\end{cases}
\]
(ii) The case $r>2$:

By Lemma~\ref{mslemma} we know that
\[
{\mathcal{L}}(E_{1},\ldots,E_{n};\mathbb{K})=\mathcal{L}_{ms(\infty;t_{1},
\ldots, t_{n})}(E_{1},\ldots,E_{n};\mathbb{K}),
\]
where we choose for $i=1, \ldots,n$
\[
\frac{1}{t_{i}}=
\begin{cases}
\frac{1}{2} & \text{if $E_{i}$ is of cotype~$2$},\\
\frac{1}{q_{i,0}} & \text{if $E_{i}$ is of cotype~$q_{i}>2$ and $1 \le
q_{i,0}<q_{i}^{\ast}$.}%
\end{cases}
\]
By Theorem~\ref{bohnenblust2}
\[
\mathcal{L}(E_{1},\ldots,E_{n};{\mathbb{K}})=\mathcal{L}_{ms(2;p_{1}^{(n)},
\ldots, p_{n}^{(n)})}(E_{1},\ldots,E_{n};{\mathbb{K}})
\]
where for $i=1, \ldots,n$
\[
\frac{1}{p_{i}^{(n)}}=
\begin{cases}
\frac{1}{2}+\frac{1}{2}-\frac{1}{2n} & \text{if $E_{i}$ is of cotype~$2$},\\
\frac{1}{q_{i,0}} +\frac{1}{{q_{i,0}}^{\ast}}-\frac{1}{{q_{i,0}}^{\ast}n} &
\text{if $E_{i}$ is of cotype~$q_{i}>2$ and $1 \le q_{i,0}<q_{i}^{\ast}$.}%
\end{cases}
\]
Now Proposition~\ref{mssandwichmixed} with $\theta=2/r$ gives the claim.
\end{proof}

\begin{remark}
\begin{enumerate}
[(a)]

\item Note that in the linear case, $\mathcal{L}_{as(q_{0};p_{0})}%
\subseteq\mathcal{L}_{as(q_{1};p_{1})}$ whenever $q_{0}\leq q_{1}$, $p_{0}\leq
p_{1}$ and $\frac{1}{p_{0}}-\frac{1}{q_{0}}\leq\frac{1}{p_{1}}-\frac{1}{q_{1}%
}$. Moreover, if $E$ is of cotype~$2$, then $\mathcal{L}_{as(q_{0};p_{0}%
)}(E;F)=\mathcal{L}_{as(q_{1};p_{1})}(E;F)$ if additionally $p_{0}\leq
p_{1}\leq2$ and $\frac{1}{p_{0}}-\frac{1}{q_{0}}=\frac{1}{p_{1}}-\frac
{1}{q_{1}}$. It is not known if there is any close analog for multiple
$(q;p)$-summing operators, but we observe that if in the above corollary all
spaces involved are of cotype~$2$, then the indices associated satisfy
$\frac{1}{r_{i}^{(n)}}-\frac{1}{r}=\frac{n-1}{2n}$ for all $\frac{2n}{n+1}\leq
r\leq2$. So, if there was some similar inclusion formula at least for multiple
$(q;p)$-summing multilinear forms, the above result for $r\leq2$ would
immediately follow from the original Bohnenblust-Hille result. In Corollary
\ref{426a} \ we will show that there is no such result for the case of
bilinear forms on Hilbert spaces.

\item For $r>2$, the above corollary shows in which sense
Theorem~\ref{bohnenblust2} improves upon Proposition~\ref{bhcotype} for the
special case $F={\mathbb{K}}$ and $q=2$.

\item A natural question is whether the estimates of Corollary
\ref{bohnenblust3} are optimal. In the next section we show that, in some
sense, for the case of Hilbert spaces our results are optimal.
\end{enumerate}
\end{remark}

We continue with a statement where spaces with finite cotype are mixed up with
arbitrary spaces.

\begin{corollary}
Let $E_{1}, \ldots, E_{n}$ be Banach spaces with finite cotype, $E_{n+1},
\ldots E_{n+k}$ be $\mathcal{L}_{\infty}$-spaces and $E_{n+k+1}, \ldots,
E_{n+k+\ell}$ be arbitrary Banach spaces, and let $\frac{2n}{n+1} \le r
<\infty$. Then
\[
\mathcal{L}(E_{1},\ldots, E_{n+k+\ell};{\mathbb{K}})=\mathcal{L}%
_{ms(r;r_{1}^{(n)}, \ldots, r_{n}^{(n)}, r-\varepsilon, \ldots, r-\varepsilon,
1, \ldots, 1)}(E_{1},\ldots,E_{n+k+\ell};{\mathbb{K}})
\]
for all $\varepsilon>0$ sufficiently small and $r_{i}^{(n)}$, $i=1, \ldots,n$
as in Corollary~\ref{bohnenblust3}.
\end{corollary}

\begin{proof}
By iteration, Corollary~\ref{bohnenblust3} and Proposition~\ref{w} give
\[
\mathcal{L}(E_{1},\ldots, E_{n+k+\ell};{\mathbb{K}})=\mathcal{L}%
_{ms(r;r_{1}^{(n)}, \ldots, r_{n}^{(n)}, 1, \ldots, 1, 1, \ldots, 1)}%
(E_{1},\ldots,E_{n+k+\ell};{\mathbb{K}}).
\]
By Lemma~\ref{cklemma}, one can improve the indices associated to the
$\mathcal{L}_{\infty}$-spaces to be as close to $r$ as wanted.
\end{proof}

\subsection{Optimality of the Bohnenblust-Hille type results}

Finally we show that the above results are partially optimal. This is
essentially based on a multilinear version of Chevet's inequality.

\begin{lemma}
\label{optimality} For all $n \in{\mathbb{N}}$ there exists a constant
$d_{n}>0$ such that for all $m \in{\mathbb{N}}$ there exists an $n$-linear
form $\varphi_{m}: l_{2}^{m} \times\cdots\times l_{2}^{m} \rightarrow
{\mathbb{K}}$ with $\|\varphi_{m}\| \le d_{n} \, m^{\frac{1}{2}}$ of the form
\[
\varphi_{m}= \sum_{j_{1}, \ldots,j_{n}=1}^{m} \varepsilon_{j_{1},\ldots,j_{n}}
e_{j_{1}} \otimes\cdots\otimes e_{j_{n}},
\]
where $\varepsilon_{j_{1},\ldots,j_{n}} \in\{-1,1\}$.
\end{lemma}

\begin{proof}
First recall the identification $\mathcal{L}(^{n}(l_{2}^{m});{\mathbb{K}%
})=l_{2}^{m}\otimes_{\varepsilon}\cdots\otimes_{\varepsilon}l_{2}^{m}$ (isometrically), where
the latter denotes the $n$-fold injective tensor product of $l_{2}^{m}$. Now
take a family $(g_{j_{1},\ldots,j_{n}})$ of independent standard Gaussian
random variables. Then by the $n$-linear version of Chevet's inequality in
\cite[Lemma~6]{ddgm01} and by \cite[p.~329]{Nicole} there exist $\tilde{d}_{n}>0$ and $\kappa>0$ such that
\[
\int\Vert\sum_{j_{1},\ldots,j_{n}=1}^{m}g_{j_{1},\ldots,j_{n}}e_{j_{1}}%
\otimes\cdots\otimes e_{j_{n}}\Vert_{l_{2}^{m}\otimes_{\varepsilon}%
\cdots\otimes_{\varepsilon}l_{2}^{m}}d\mu\leq\tilde{d}_{n}\int\Vert\sum
_{i=1}^{m}g_{i}e_{i}\Vert_{l_{2}^{m}}d\mu{\leq}\kappa\tilde{d}_{n}m^{\frac{1}{2}}.
\]
Now take a family $(\varepsilon
_{j_{1},\ldots,j_{n}})$ of independent Bernouilli random variables. It is
well-known that (up to a universal constant) Bernouilli averages are dominated
by Gaussian averages, thus there exists a constant $d_{n}>0$ such that
\[
\int\Vert\sum_{j_{1},\ldots,j_{n}=1}^{m}\varepsilon_{j_{1},\ldots,j_{n}%
}e_{j_{1}}\otimes\cdots\otimes e_{j_{n}}\Vert_{l_{2}^{m}\otimes_{\varepsilon
}\cdots\otimes_{\varepsilon}l_{2}^{m}}d\mu\leq{d}_{n}m^{\frac{1}{2}}.
\]
Now since these Bernouilli averages are dominated (up to the constant $d_{n}$)
by $m^{\frac{1}{2}}$, there exists an $n$-linear form $\varphi_{m}$ as desired.
\end{proof}

This now gives the following partially optimal statement for $n$-linear forms
on Hilbert spaces:

\begin{theorem}
\label{hilbert}
Let $H_{1},\ldots,H_{n}$ be Hilbert spaces and $\frac{2n}{n+1}\leq r<\infty$.
Then
\[
\mathcal{L}(H_{1},\ldots,H_{n};{\mathbb{K}})=\mathcal{L}_{ms(r;r_{n})}%
(H_{1},\ldots,H_{n};{\mathbb{K}}),
\]
where
\[
\frac{1}{r_{n}}=\frac{1}{2}+\frac{1}{r}-\frac{1}{\max(r,2)n}.
\]
For $\frac{2n}{n+1}\leq r\leq2$, the parameter $r_{n}$ is best possible.
\end{theorem}

\begin{proof}
By Corollary~\ref{bohnenblust3} and the trivial fact that any Hilbert space is of cotype~2, it remains to show optimality for $\frac{2n}{n+1}\leq r\leq2$. Without loss of
generality we may assume that all Hilbert spaces involved are
infinite-dimensional. Assume that
\[
\mathcal{L}(H_{1},\ldots,H_{n};\mathbb{K})=\mathcal{L}_{ms(r;p_{n})}%
(H_{1},\ldots,H_{n};\mathbb{K})
\]
with $1\leq p_{n}\leq 2$. Then there exists $C_{n}\geq0$ independent of $m$
and $T\in\mathcal{L}(H_{1},\ldots,H_{n};\mathbb{K})$ such that
\begin{equation}
\left(  \sum\limits_{j_{1},...,j_{n}=1}^{m}|T(x_{j_{1}}^{(1)},...,x_{j_{n}%
}^{(n)})|^{r}\right)  ^{1/r}\leq C_{n}\,\Vert T\Vert\,%
{\textstyle\prod\limits_{k=1}^{n}}
\left\Vert (x_{j}^{(k)})_{j=1}^{m}\right\Vert _{w,p_{n}},\label{dummy}%
\end{equation}
for every $m\in\mathbb{N}$ and $x_{j}^{(k)}\in H_{k}$, $j=1,...,m$ and
$k=1,...,n.$ Now we may assume that $H_{k}=l_{2}^{m}$, and $x_{j}^{k}=e_{j}$.
Then the right-handside in \eqref{dummy} equals
\[
C_{n}\,\Vert T\Vert\,\Vert\text{id}:l_{2}^{m}\rightarrow l_{p_{n}}^{m}%
\Vert^{n}=C_{n}\Vert T\Vert\,m^{\frac{n}{p_{n}}-\frac{n}{2}}.
\]
Now by Lemma~\ref{optimality} there exists an $n$-linear form $\varphi
_{m}:l_{2}^{m}\times\cdots\times l_{2}^{m}\rightarrow{\mathbb{K}}$ with
$\Vert\varphi_{m}\Vert\leq d_{n}\,m^{\frac{1}{2}}$ of the form
\[
\varphi_{m}=\sum_{j_{1},\ldots,j_{n}=1}^{m}\varepsilon_{j_{1},\ldots,j_{n}%
}e_{j_{1}}\otimes\cdots\otimes e_{j_{n}},
\]
where $\varepsilon_{j_{1},\ldots,j_{n}}\in\{-1,1\}$. For this $n$-linear form,
the right-handside in \eqref{dummy} now can be estimated from above by
$\tilde{C}_{n}\,m^{\frac{n}{p_{n}}-\frac{n}{2}+\frac{1}{2}}$ for some other
constant $\tilde{C}_{n}>0$ independent of $m$ and $\varphi_{m}$, whereas the
left-handside equals $m^{\frac{n}{r}}$. Thus,
\[
\frac{1}{p_{n}}\geq\frac{1}{2}+\frac{1}{r}-\frac{1}{2n}=\frac{1}{r_{n}}%
\]
and hence $p_{n}\leq r_{n}$.
\end{proof}

We conjecture that our result in the case $r > 2$ is also optimal; this
conjecture is motivated by the upcoming proposition which deals with the case
$n=2$. 

\begin{proposition}
\label{asoptimality}
Let $2 \le r <\infty$. Then for $1 \le p \le r$ the following are equivalent:

\begin{enumerate}
[(i)]

\item $\mathcal{L}(l_{2}, l_{2};{\mathbb{K}})=\mathcal{L}_{ms(r;p)}(l_{2},
l_{2};{\mathbb{K}})$;

\item $\mathcal{L}(l_{2}, l_{2};{\mathbb{K}})=\mathcal{L}_{as(r;p)}(l_{2},
l_{2};{\mathbb{K}})$;

\item $\frac{1}{p} \ge\frac{1}{2}+\frac{1}{2r}$.
\end{enumerate}
\end{proposition}

\begin{proof}
By Corollary~\ref{bohnenblust3} and the inclusion $\mathcal{L}_{ms(r;p)}
\subseteq\mathcal{L}_{as(r;p)}$ we only have to show that (ii) implies (iii).
Assume that (ii) holds, and consider the bilinear operator $T(x,y):=
\sum_{i=1}^{\infty}x_{i} y_{i}$, $x,y \in l_{2}$, which is of norm $1$ by the
H\"older inequality. Then by definition of the absolutely
$(r;p)$-summing norm and (ii) there exists $C>0$ (independent of $m$) such
that
\[
\left(  \sum\limits_{j=1}^{m}|T(x_{j}^{(1)},x_{j}^{(2)})|^{r}\right)
^{1/r}\leq C \,%
\left\Vert (x_{j}^{(1)})_{j=1}^{m}\right\Vert _{w,p} \left\Vert (x_{j}%
^{(2)})_{j=1}^{m}\right\Vert _{w,p} ,
\]
for every $m\in\mathbb{N}$ and $x_{j}^{(k)}\in E_{k}$, $j=1,...,m$ and
$k=1,2.$ Now choose $x_{j}^{(k)}=e_{j}$ for $j=1,...,m$ and $k=1,2.$ Then
\[
\left( \sum_{j=1}^{m} |T(x_{j}^{(1)},x_{j}^{(2)})|^{r}\right) ^{1/r}=m^{1/r}
\]
and
\[
\left\Vert (x_{j}^{(1)})_{j=1}^{m}\right\Vert _{w,p} \left\Vert (x_{j}%
^{(2)})_{j=1}^{m}\right\Vert _{w,p}=\|\text{id}: l_{2}^{m} \rightarrow
l_{p}^{m}\|^{2}=m^{\max(\frac{2}{p}-1,0)}.
\]
This now implies $p \le 2$ and $\frac{1}{r} \le\frac{2}{p} -1 $, and therefore (iii).
\end{proof}

With regard to the remark after Corollary~\ref{bohnenblust3}, this immediately
implies the following supplement to the non-existence of a general inclusion
result of the type $\mathcal{L}_{ms,q_{0}} \subseteq\mathcal{L}_{ms,q_{1}}$
whenever $1 \le q_{0} < q_{1} <\infty$ in \cite{inclusion} -- note that
nevertheless $\mathcal{L}_{ms,1}(H_{1}, \ldots, H_{n};{\mathbb{K}%
})=\mathcal{L}_{ms,q}(H_{1}, \ldots, H_{n};{\mathbb{K}})$ for all $1
<q<\infty$ by \cite{inclusion} and hence, $\mathcal{L}_{ms,q_{0}}(H_{1},
\ldots, H_{n};{\mathbb{K}}) \subseteq\mathcal{L}_{ms,q_{1}}(H_{1}, \ldots,
H_{n};{\mathbb{K}})$ whenever $1 \le q_{0} < q_{1} <\infty$.

\begin{corollary}
\label{426a}There does not exist a general inclusion result of the type
$\mathcal{L}_{ms(q_{0};p_{0})}\subseteq\mathcal{L}_{ms(q_{1};p_{1})}$ whenever
$p_{0}<p_{1}$ and $\frac{1}{p_{0}}-\frac{1}{q_{0}}\leq\frac{1}{p_{1}}-\frac
{1}{q_{1}}$, not even for bilinear forms on Hilbert spaces.
\end{corollary}

\begin{proof}
Assume that such a general inclusion result would hold. Then the original
Bohnenblust-Hille result would imply that $\mathcal{L}(l_{2},l_{2}%
;{\mathbb{K}})=\mathcal{L}_{ms(r;p)}(l_{2},l_{2};{\mathbb{K}})$ whenever $r>2$
and $1\leq p\leq r$ are such that
\[
1-\frac{1}{4/3}\leq\frac{1}{p}-\frac{1}{r}\text{.}%
\]
However, since $r>2$, this would contradict $\frac{1}{p}\geq\frac{1}{2}%
+\frac{1}{2r}$ from the above proposition (for example, when $r=4$ and $p=2$).
\end{proof}

Optimality of coincidences for multiple summing operators on Hilbert spaces has a great effect on optimality of coincidences for operators on so-called $K$-convex spaces. For this notion, we refer to \cite{DJT}; a fundamental characterization due to G. Pisier says that a Banach space is $K$-convex if and only if it is of non-trvial type (see. e.g., \cite[Theorem~13.3]{DJT}).   

\begin{lemma}
\label{banachoptimal}
Let $E_1, \ldots, E_n$ be $K$-convex Banach spaces, $F$ be an arbitrary Banach space, $1 \le r <\infty$ and $1 \le r_1, \ldots,r_n \le r$. Then $\mathcal{L}(E_1, \ldots,E_n;F)=\mathcal{L}_{ms(r;r_1, \ldots,r_n)}(E_1, \ldots,E_n;F)$ implies $\mathcal{L}(^n l_2;F)=\mathcal{L}_{ms(r;r_1, \ldots,r_n)}(^n l_2;F)$.
\end{lemma}
\begin{proof}
This follows by standard arguments from the fact that a Banach space is $K$-convex if and only if it contains the $l_2^n$'s uniformly and uniformly complemented (see, e.g., \cite[Theorem~19.3]{DJT}).
\end{proof}

We now can state in which sense some of our Bohnenblust-Hille type theorems are optimal also for multilinear forms on arbitrary Banach spaces; note that, e.g., $l_p$ for $1<p\le 2$ is of cotype~$2$ as well as $K$-convex. 

\begin{corollary}
For $\frac{2n}{n+1} \le r \le 2$ and $E_1, \ldots, E_n$ all $K$-convex Banach spaces of cotype~$2$, the result from Corollary~\ref{bohnenblust3} is best possible.
\end{corollary}
\begin{proof}
This immediately follows from Theorem~\ref{hilbert} and Lemma~\ref{banachoptimal}.
\end{proof}

\begin{remark}
The space $l_1$ is of cotype~$2$ but not $K$-convex, and Corollary~\ref{bohnenblust3} is far from being best possible for multilinear forms on $^n l_1$: Every multilinear operator from $^n l_1$ into a Hilbert space is multiple $r$-summing for all $1 \le r \le 2$ (see \cite{bpgv}). Thus, the $K$-convexity condition for optimality is not superfluous and seems to be quite appropriate.  
\end{remark}

We finish with the following generalization of \cite[Remark 2.1]{ckp}, where the case $r=2$ is treated. Although it can be proved by other means, it shows in which sense our more abstract coincidence result for multiple summing operators can be used to obtain inequalities more closely related to the original Bohnenblust-Hille inequality.  

\begin{corollary}
Let $T \in \mathcal{L}(^n l_2;\K)$ and $1 \le r <\infty$. Then for some constant $C_n>0$ not depending on $m$ and $T$, the following hold:
\begin{enumerate}[(i)]
\item
If $1 \le r \le 2$, then
\[
\left(  \sum\limits_{j_{1},...,j_{n}=1}^{m}|T(e_{j_{1}},...,e_{j_{n}%
})|^{r}\right)  ^{1/r}\leq C_n \, m^{\frac{n}{r}-\frac{1}{2}} \norm{T},  
\]
and the exponent $\frac{n}{r}-\frac{1}{2}$ is best possible.
\item
If $2 \le r <\infty$, then
\[
\left(  \sum\limits_{j_{1},...,j_{n}=1}^{m}|T(e_{j_{1}},...,e_{j_{n}%
})|^{r}\right)  ^{1/r}\leq C_n \, m^{\frac{n-1}{r}} \norm{T},  
\]
and the exponent $\frac{n-1}{r}$ is best possible when $n=2$.
\end{enumerate}
\end{corollary}
\begin{proof}
(ii) The estimate follows from Theorem~\ref{hilbert} and similar reasoning as in its proof. Alternatively, one may use \cite[Remark 2.1]{ckp} (the case $r=2$) and the fact that for $2<r<\infty$ the space $l_r^{m^n}$ is of power type $1-\frac{2}{r}$ with respect to the interpolation couple $(l_2^{m^n},l_\infty^{m^n})$; note that this even shows that one may choose $C_n=1$. The optimality for $n=2$ can be seen using the same bilinear form as in the proof of Proposition~\ref{asoptimality}.

(i) The estimate follows by factorization through $l_2^{m^n}$ from the case $r=2$. The optimality can be seen by using the $n$-linear form from Lemma~\ref{optimality}.
\end{proof}

\appendix

\section{Complexification}

The following definitions and results are essentially based on ideas presented
in \cite[p.~68--70]{David}, and we omit the mostly straightforward proofs of
the results.

Let $E$ be a real Banach space, and define the complex vector space $\tilde
{E}=E \oplus E$ with the operations
\begin{align*}
(x,y)+(u,v)  &  =(x+u,y+v), \qquad x,y,u,v \in E,\\
(\alpha+i\beta)(x,y)  &  =(\alpha x-\beta y,\beta x+\alpha y), \qquad x,y \in
E, \quad\alpha,\beta\in{\mathbb{R}}.
\end{align*}
This becomes a complex Banach space under the norm
\[
\|x+iy\|_{\tilde{E}}=\|x \otimes e_{1} +y\otimes e_{2}\|_{E \otimes_{\pi}%
l_{2}^{2}}, \qquad x+iy=(x,y) \in\tilde{E},
\]
where $E \otimes_{\pi}l_{2}^{2}$ denotes the projective tensor product of $E$
with $l_{2}^{2}$.

If $T \in\mathcal{L}(E_{1}, \ldots,E_{n};F)$ is an operator between real
Banach spaces, we define its complexification $\tilde{T} \in\mathcal{L}%
(\tilde{E}_{1}, \ldots, \tilde{E}_{n}; \tilde{F})$ by
\[
\tilde{T}(x^{1,0}+ix^{1,1}, \ldots, x^{n,0}+ix^{n,1})= \sum_{\varepsilon_{1},
\ldots, \varepsilon_{n}} i^{\sum_{k=1}^{n} \varepsilon_{k}} T(x^{1,\varepsilon
_{1}}, \ldots, x^{n,\varepsilon_{n}}).
\]

With these definitions, one can easily prove the following:

\begin{proposition}
\begin{enumerate}
[(a)]

\item Let $T \in\mathcal{L}(E_{1}, \ldots,E_{n};F)$ be an operator between
real Banach spaces and $\tilde{T}$ its complexification. If $1 \le p_{1},
\ldots, p_{n},q \le\infty$ are such that $1/q \le1/{p_{1}} + \ldots1/{p_{n}}$,
then $T$ is absolutely $(q;p_{1},\ldots,p_{n})$-summing if and only if
$\tilde{T}$ is.

\item Let $T \in\mathcal{L}(E_{1}, \ldots,E_{n};F)$ be an operator between
real Banach spaces and $\tilde{T}$ its complexification. If $1 \le p_{1},
\ldots, p_{n} \le q \le\infty$, then $T$ is multiple $(q;p_{1},\ldots,p_{n}%
)$-summing if and only if $\tilde{T}$ is.
\end{enumerate}
\end{proposition}

The following is only a short list of properties of a Banach space which are
stable under complexification, essentially the ones that we need for our purposes.

\begin{proposition}
The following properties of a Banach space are stable under complexification:

\begin{enumerate}
[(a)]

\item having cotype $q$, for $2 \le q <\infty$;

\item being an $\mathcal{L}_{p}$-space, for $1 \le p \le\infty$.
\end{enumerate}
\end{proposition}



\begin{thebibliography}{99}                                                                                               %
\bibitem{alencarmatos} R. Alencar and M.C. Matos, \emph{Some classes of multilinear mappings between Banach spaces}, 1989, Publicaciones del Departamento de An\'alisis Mathem\'atico, sec. 1, no. 12. Universidad Complutense de Madrid.

\bibitem {Bergh}J. Bergh and J. L\"{o}fstr\"{o}m, \emph{Interpolation spaces},
Springer-Verlag, 1976.

\bibitem {BBPR}O. Blasco, G. Botelho, D. Pellegrino and P. Rueda,
\emph{Summability of multilinear mappings: Littlewood, Orlicz and beyond}, preprint.

\bibitem {BH}H. F. Bohnenblust and E. Hille, \emph{On the absolute convergence
of Dirichlet series}, Ann. Math. \textbf{32} (1931), 600--622.

\bibitem {bpgv}F. Bombal, D. Per\'{e}z-Garc\'{\i}a, and I. Villanueva,
\emph{Multilinear extensions of Grothendieck's theorem}, Quart. J. Math.
\textbf{55} (2004), 441--450.

\bibitem {B}G. Botelho, \emph{Cotype and absolutely summing multilinear
mappings and homogeneous polynomials}, Proc. Royal Irish Acad. 97 (1997), 145--153.

\bibitem {BBJP-PAMS}G. Botelho, H.-A. Braunss, H. Junek, and D.
Pellegrino,\emph{ Inclusions and coincidences for multiple summing multilinear
mappings}. Proc. Amer. Math. Soc., to appear.

\bibitem {bparxiv1}G. Botelho and D. Pellegrino, \emph{When every multilinear
mapping is multiple summing}, Math. Nachr., to appear.

\bibitem {bparxiv2}G. Botelho and D. Pellegrino, \emph{Coincidences for
multiple summing mappings}, arXiv:0809.4171v2.

\bibitem{ckp} F. Cobos, T. K\"uhn, and J. Peetre, \emph{On $\mathfrak{S}_p$-classes of trilinear forms}, J. London Math. Soc. \textbf{59} (1999), 1003--1022.

\bibitem {ddgm01}A. Defant, D\'{\i}az, D. Garc\'{\i}a, and M. Maestre,
\emph{Unconditional basis and Gordon-Lewis constants for spaces of
polynomials}, J. Funct. Anal. \textbf{181} (2001), 119--145.

\bibitem {DM}A. Defant and C. Michels, \emph{A complex interpolation formula
for tensor products of vector-valued Banach function spaces}, Arch. Math.
\textbf{74} (2000), 441--451.

\bibitem {DS}A. Defant and P. Sevilla-Peris, \emph{A new multilinear insight
on Littlewood's $4/3$-inequality}, J. Funct. Anal., to appear.

\bibitem {DJT}J. Diestel, H. Jarchow and A. Tonge, \emph{Absolutely summing
operators}, Cambridge University Press 1995.

\bibitem {JMP}H. Junek, M. Matos and D. Pellegrino, \emph{Inclusion theorems
for absolutely summing holomorphic mappings}, Proc. Amer. Math. Soc.
\textbf{136} (2008), 3983--3991.

\bibitem{littlewood} J.E. Littlewood, \emph{On bounded bilinear forms in an infinite number of variables}, Quart. J. Math. \textbf{1} (1930), 164--174.

\bibitem{matos} M.C. Matos, \emph{On multilinear mappings of nuclear type}, Rev. Mat. Univ. Complut. Madrid \textbf{6} (1993), 61--81.

\bibitem {Michels-Belg}C. Michels, \emph{One-sided interpolation of injective
tensor products of Banach spaces}, Bull. Belg. Math. Soc. Simon Stevin
\textbf{14} (2007), 531--538.

\bibitem {Pell-Irish}D. Pellegrino, \emph{Cotype and nonlinear absolutely summing
mappings}, Proc. Roy. Irish Acad. \textbf{105A} (2005), 75--91.

\bibitem {David}D. P\'{e}rez-Garc\'{\i}a, \emph{Operadores multilineales
absolutamente sumantes}, Doctoral Thesis, Universidad Complutense de Madrid, 2003.

\bibitem {inclusion}D. P\'{e}rez-Garc\'{\i}a, \emph{The inclusion theorem for
multiple summing operators}, Studia Math. \textbf{165} (2004), 275--290.

\bibitem {PPGG}D. P\'{e}rez-Garc\'{\i}a, \emph{The trace class is a
$Q$-algebra}, Ann. Acad. Sci. Fenn. Math. \textbf{21} (2006), 287--295.

\bibitem {pgvck}D. Per\'{e}z-Garc\'{\i}a and I. Villanueva, \emph{Multiple
summing operators on $C(K)$-spaces}, Ark. Mat.~\textbf{42} (2004), 153--171.

\bibitem {popa}D. Popa, \emph{Reverse inclusions for multiple summing
mappings}, J. Math. Anal. Appl., to appear.

\bibitem {souza}M. L. V. Souza, \emph{Aplica\c{c}\~{o}es multilineares
completamente absolutamente somantes}, Doctoral Thesis, Unicamp, 2003.

\bibitem {Nicole}N. Tomczak-Jaegermann, Banach-Mazur distances and
finite-dimensional operators ideals, Longman Scientific \& Technical,
Harlow/Londin 1989.
\end{thebibliography}
\end{document}